# Optimal Trajectory-Planning of UAVs via B-Splines and Disjunctive Programming

Alireza Babaei and Amirhossein Karimi

*Abstract*—This paper investigates an efficient algorithm for trajectory planning problem of autonomous unmanned aerial vehicles which fly over three-dimensional terrains. The proposed algorithm combines convex optimization with disjunctive programming and receding horizon concept, which has many advantages, such as a high computational speed. Disjunctive programming is applied in order to relax the non-convex constraints of the problem. Moreover, the B-spline curves are employed to represent the trajectories which should be generated in the optimization process.

## I. INTRODUCTION

UNMANNED Aerial Vehicles (UAVs) [1] have become increasingly attractive due to the unique applications of them in environments which are mostly inaccessible to other types of aircrafts, as in the case of flying at low altitude for the purpose of urban operations and the missions in which human presence is hazardous and thus impossible [2]. Since the UAVs operate without onboard pilots, there are two types of controlling strategy for them, i.e., remotely or autonomously. In this paper the second is intended in which an online planner is required that generates the optimal trajectories for a UAV in every potential situation.

Recent advances in the navigation systems, microcontroller design, digital computers, optimization techniques, and digital vision systems have broadened the versatility and autonomy of the UAVs with a high degree of reliability in their operations [3]. Various applications of UAVs in civil, military, and commercial ones incorporate surveillance, weather and atmospheric monitoring, reconnaissance, emergency communications, environmental and meteorological monitoring, telecommunications, aerial photography, border patrol, search-and-rescue tasks, etc. [4].

Among several open problems in the field of autonomous UAVs, trajectory generation is of immense significance, in which an optimal or near-optimal feasible flight trajectory should be generated automatically. A trajectory is feasible if it respects all the limitations and constraints imposed by the physical characteristics of the vehicle and the environment (such as minimum turning radius, minimum and maximum speed, and avoiding the obstacles and terrain collision) and mission constraints (such as passing through the specified waypoints). There are a variety of studies, in which different constraints and objective functions for the UAV trajectory generation problem have been taken into consideration, such as implementing special tasks [5], obstacle avoidance [6], curvature constraints [3], multiple UAVs cooperation [7], etc.

Generally, the optimization problem of trajectory generation is NP-hard (Non-deterministic Polynomial-time hard), which is widely regarded as a sign that a polynomial-time algorithm for this problem is unlikely to exist [8]. Therefore, in the previous studies, several algorithms have been contributed based on the approximate algorithms, i.e., heuristic and meta-heuristic ones [9], [10]. In heuristic approaches, the problem is usually relaxed to a simpler one, which can be solved using different methods such as mixed-integer linear programming [11] or dynamic programming [12]. Furthermore, meta-heuristic approaches such as evolutionary algorithms are widely used to solve the optimal trajectory generation problem for the UAVs [13], [14]. The main problem with evolutionary algorithms is requiring high computational burden, and therefore the evolutionary-based planners work mostly in offline mode. Moreover, if it is required for an evolutionary-based planner to work in online mode due to the presence of unexpected pop-ups (unknown threats), a simpler model of the problem and thus less optimization indexes are considered in order to speed up the algorithm [9].

Safe maneuver of the UAVs considering the terrain collision avoidance consists in two different problems, i.e., the Terrain Following (TF) and Terrain Avoidance (TA) [15]–[17]. In the TF maneuver, the UAV follows the terrain by climbing and descending in a vertical plane through the flight path vector in order to maintain as closely as possible to a desired height above the terrain, for which only the longitudinal dynamics of the UAV is engaged. However, the TA maneuver is accomplished in a horizontal plane parallel to the ground, in which the UAV flies at the low altitudes without colliding with the terrains, to the extent that in some studies [18] it is assumed that any change to the flying altitude is mathematically negligible. In the latter maneuver, lateral dynamics of the flying vehicle plays a prominent role in comparison to the longitudinal one. It should be noted that in the most real cases a synergy-based implementation of the two aforementioned maneuvers is required in practice, in which the UAV is capable of flying over the mountains as well as turning the obstacles in low altitudes allowing an improvement in the terrain masking. Therefore, in the present study, the combined Terrain Following-Terrain Avoidance (TF/TA) maneuver is taken into account for the problem of optimal trajectory planning for the UAVs.

First study on the B-spline (Basis spline) curves was conducted in [19], and then pursued by introducing developing algorithms in [20]–[22]. In the computer-aided design and computer graphics, B-splines have many useful applications such as those proposed in [23]–[25].

In this paper, an algorithm for optimal trajectory planning for UAVs will be proposed using a combination of discrete optimization, disjunctive programming, and receding horizon. Furthermore, different objective functions and constraints are considered in the problem in order to make it appropriate for various realistic scenarios. Moreover, B-spline curves are applied to represent the trajectories which should be generated in the optimization process. The computational time for the approach presented in this paper is low, which opens an avenue to use it as a systematic algorithm for online purposes. In this research, the UAV is approximated as a point mass, moving with limited speed, and thus its dynamics is expressed as a simple point mass. It should be noted that the point mass model captures the substantial features of the UAV real dynamics such as velocity and curvature constraints [26]. Additionally, it is assumed that the UAV is equipped with some sensors in order to detect unknown obstacles and terrain. In other words, the vehicle finds a new obstacle or terrain when they come within its field of detection, which is considered as the internal region of a spherical pentahedron with the range $R$ and angles $\theta_1$ and $\theta_2$ as illustrated in Fig. I.

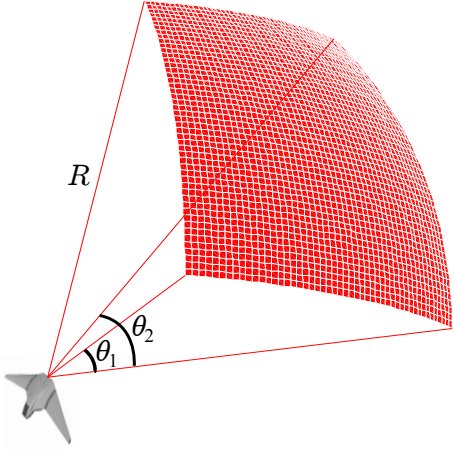

Fig. 1. Field of detection of the UAV sensors which is a spherical pentahedron.

The rest of this paper is organized as follows. Section 2 introduces the B-spline curves and their important properties. In section 3, the principles for discrete optimization and receding horizon concept have been touched briefly upon. Section 4 describes disjunctive programming and the procedure by which one can reduce it to a mixed integer linear programming. Thereafter, in section 5, the problem formulation including the objective function and the required constraints are proposed. Then, in section 6, various test scenarios have been taken into consideration and the proposed algorithm has been implemented, for which the results are illustrated. Finally, the paper concludes by providing some remarks and ongoing work.

## II. B-SPLINE CURVES

As mentioned previously, in this study B-spline curves are applied to represent the trajectories which should be generated in the optimization process. Before delving into the details of B-splines, some terms are introduced in the following [27]:

- Data point: the points that the curve should pass through them.
- Control point: the points for controlling the shape of the curve. These points allow the user to modify the shape of the curve by changing them.
- Geometric continuity: a curve has geometric continuity $G^n$ at a connecting point if every pair of the first $n$ derivatives of the two segments have identical direction at that point.
- Parametric continuity: If as well as the same direction, the derivatives have equal magnitudes at the point, then the curve have $C^n$ parametric continuity at that point.

Unlike linear and polynomial interpolations [27], one of the positive points of B-splines is high interactivity. It means that there are controllable parameters, i.e., control data, which can modify the shape of the curve rationally by the user. It should be noted that in case the order of the B-spline curve is equal to the number of control points, it reduces to a Bézier curve. There are two fundamental drawbacks in Bézier representation for which B-splines can be considered as a remedy. Firstly, it is only globally controllable, i.e., changing any control point modifies the entire curve. Secondly, the degree of a Bézier function depends on the number of control points and therefore, in order to increase the complexity of the shape of the curve by adding control points, one should increase the degree of the curve as well or apply multi-segments Bézier curves which satisfy the continuity conditions between consecutive segments [28].

Given a non-descending sequence of breaking points, i.e., knot vector $t = [t_0, ..., t_m]^T$, the B-spline basis function of degree $k-1$, defined on $[t_i, t_{i+k}]$, is obtained using the recursive Cox-DeBoor formula [27]:

$$N_{i,k}(t) = \frac{t - t_i}{t_{i+k-1} - t_i} N_{i,k-1}(t) + \frac{t_{i+k} - t}{t_{i+k} - t_{i+1}} N_{i+1,k-1}(t),$$
$$i = 0, ..., n \quad (1)$$

which starts with,

$$N_{i,1}(t) = \begin{cases} 1 & \text{if } t_i \leq t < t_{i+1} \\ 0 & \text{otherwise} \end{cases} \quad (2)$$

It should be noted that the basis functions $N_{i,k}(t)$, are strictly positive on $(t_i, t_{i+k})$ and zero otherwise. Moreover, they have the following interesting property, namely partition of unity [28]:

$$\sum_{i=0}^{n} N_{i,k}(t) = 1, \ \forall t \in [0, t_m] \quad (3)$$

A B-spline curve ($\mathbf{r}(t)$) of order $k$ can be expressed as a linear combination of B-spline basis functions using the control points ($\mathbf{q}_i$) as coefficients:

$$\mathbf{r}(t) = \sum_{i=0}^{n} \mathbf{q}_i N_{i,k}(t) \quad (4)$$

By normalizing the knot vector, one can enforce the B-spline curve to cover the interval [0,1], which is more appropriate for the purpose of higher numerical accuracy [29].

In the following some properties of B-spline curves have been introduced:

- B-spline curves do not in general pass through the control points, but to define the shape of the curve.
- At the endpoints of the interval $[0, t_m]$, the B-spline function $\mathbf{r}(t)$ has zero values. Therefore, if a B-spline curve defined on $[0, t_m]$ is desired to have non-zero values at the endpoints, it is necessary to place $k$ knots at each location, i.e., increasing the multiplicity of the knots at the endpoints.
- If $\mathbf{r}(t)$ is a B-spline of degree $k-1$ without multiplicity, i.e., the knot vector is strictly increasing, then $\mathbf{r}(t)$ and its derivatives up to order $k-2$ are continuous on $[0, t_m]$ [30].
- Increasing the multiplicity of a knot is equivalent with the reduction of the continuity of the curve at that knot. In other words, at a knot with multiplicity $l$ ($l \leq k$), only the derivatives of $\mathbf{r}(t)$ up to order $k-l-1$ are continuous, i.e., $C^{k-l-1}$ continuity [28].
- The number of control points $(n+1)$ plus the order of the curve $(k)$ is equal with the number of elements of knot vector, namely $m = n + k$.
- Convex hull property: The B-spline curve lies within the convex hull of the control points [28]. The convex hull of a set of points is the smallest convex set containing all those points [31]. This property stems from the partition of unity.

Finally, It should be noted that since discontinuity in the acceleration vector is undesirable, the generated motion trajectory should be at least $C^2$ continuous. Therefore, the B-splines used in this paper will be constrained to third-degree B-splines to ensure the $C^2$ continuity of the trajectory while require less computation complexity than that of using higher degree curves.

## III. Discrete Optimization and Receding Horizon

Using the concept of discrete programming and receding horizon can greatly reduce the the computational time and leads to an appropriate algorithm for the online applications [32]. Discrete optimization means taking the optimal trajectory as a series of discrete points with equal or unequal time intervals. In this way, the number of variables is very large and thus increases the computational time. To solve this problem, the concept of receding horizon [33] is employed. This means that at each stage, the optimal trajectory for a temporal horizon (or spatial horizon) is obtained, but only the values obtained for the first step are retained, called state in this paper, and the rest are eliminated. This continues until the UAV reaches the target. This will speed up the calculations and, due to its predictive nature, it can be used in unknown environments and in case of environmental changes. In this study, the detectable range for each UAV is assumed to be larger than the planning horizon in order to force the trajectory to be generated according to the known data.

One of the most important topics in this field is safe receding horizon. Indeed, this expresses that the feasibility of the optimization problem over $h$ time steps started from the time step $i$, does not necessarily ensure the feasibility over $h$ time steps started from $i+1$ [33]. This may result in a collision with the obstacle at the time interval between time step $i+h$ and $i+h+1$. This led to the concept of rescue path in some studies, in which the vehicle is allowed to move to the next state, only if there exists a rescue path starting from that state [33]. It should be noted that this issue is more important for those vehicles which move at high speed, such as aircraft. In this study, we assume that the time required for the UAV to move to its safe state is sufficiently low in order to avoid collision. Furthermore, the detectable range for each UAV is assumed to be larger than the planning horizon.

In the rest of this paper, the position vector of the UAV is represented via a three-dimensional B-spline curve as follows:

$$\mathbf{x}(t) = \mathbf{r}(t) = \sum_{l=0}^{n} \mathbf{q}_l N_{l,k}(t) \tag{5}$$

where the control points, $\mathbf{q}_l (l = 0, ..., n)$, are considered as the optimization variables, and the knot vector is assumed as,

$$t = [\; \underbrace{t_0, ..., t_{k-1}}_{k \text{ equal knots}}, \; \underbrace{t_k, ..., t_n}_{(n-k+1) \text{ internal knots}}, \underbrace{t_{n+1}, ..., t_{n+k}}_{k \text{ equal knots}} ]^T \tag{6}$$

with multiplicity at the endpoints as stated previously. Furthermore, the symbol $\mathbf{x}_i = [x_i, y_i, z_i]^T$ indicates the position vector of the UAV at the $i$th time step. Moreover, the time horizon is represented with $t_h$ and the number of time steps considered for the receding horizon is shown by $h$. Discrete time interval is also displayed with $T_s$. It should be noted that each time step is equivalent with a knot in the knot vector.

## IV. Disjunctive Programming

A disjunctive programming problem required in this study can be defined as follows [34]:

$$\begin{aligned} \text{Min.} \quad & f(\mathbf{x}) \\ & \bigvee_{i=1,...,m} C_i(\mathbf{x}) \leq 0 \end{aligned} \tag{7}$$

where $\mathbf{x}$ is the vector of decision variables and $f(\mathbf{x})$ is the cost function which should be minimized. Moreover, the constraint is a disjunction of $m$ inequalities $C_i(\mathbf{x})$, which are affine functions of $\mathbf{x}$ in this paper, namely,

$$C_i(\mathbf{x}) = \mathbf{a}_i^T \mathbf{x} - c_i \tag{8}$$

where $\mathbf{a}_i$ and $c_i$ are the vector and scalar coefficients, respectively. The general definition for disjunctive programming can be found in [35]. This problem can be solved by reformulating it as a Mixed Integer Program (MIP), which can be done in different ways, including big-M method and convex hull reformulations [36], [37]. In this research, the big-M method is applied to modify the related disjunctive problem to an MIP via reformulating the aforementioned constraint as,

$$\begin{aligned} & C_i(\mathbf{x}) \leq (1 - b_i) M, \; i = 1, ..., m \\ & \sum_{i=1}^{m} b_i \geq 1 \end{aligned} \tag{9}$$

where $M$ is a sufficiently large constant and $b_i$ $(i = 1, ..., m)$ are binary variables. It is noteworthy that the second constraint in Eq. (9) is readily reduced to:

$$\sum_{i=1}^{m} b_i = 1 \qquad (10)$$

This problem can be solved efficiently using available solvers such as CPLEX [38] or Gurobi [39].

## V. PROBLEM FORMULATION

### A. Objective Function

Different terms, including travel time, path length, flight altitude, energy consumption, or a combination of them can be included in the objective function. It should be noted that in this study it is focused on the convex objective functions. In the rest of this paper, the objective function is considered as a combination of the following terms:

*1) Minimum trajectory length:* Minimizing the length of the trajectory traveled by the UAV has the advantages of diminishing the flight duration and the chance to involve with an unknown threat [9]. Furthermore, It can reduce the fuel usage, however, it is not a complete reasoning because there are other factors which affect the fuel consumption such as aerodynamics, flight performance characteristics, and environmental effects like prevailing winds [40].

To achieve this objective, during any horizon while the endpoint should approach the target state, the line segments passing through the midpoints should have the overall minimum length. Assuming the starting and final states of the UAV as $\mathbf{x}_s$ and $\mathbf{x}_f$, respectively, the objective function ($J$) can be expressed as follows:

$$J_1 = w_1 \sum_{i=0}^{h-1} ||\mathbf{x}_{i+1} - \mathbf{x}_i||_2 + w_2||\mathbf{x}_h - \mathbf{x}_f||_2 \qquad (11)$$

where $w_1$ and $w_2$ are the weights assigned to each term.

*2) Minimum flight altitude:* Flying at low altitude has the benefit of not be detected by unknown radars despite increasing the fuel consumption. To do this end, the accumulated altitude of the UAV during the time horizon is minimized, namely,

$$J_2 = w_3 \sum_{i=1}^{h} z_i \qquad (12)$$

where $w_3$ is the weight assigned to this term.

*3) Minimum number of active knots:* The $k$th-order derivative of a B-spline curve is zero throughout the whole time interval, except at the knots which are active, namely the knots for which the $k$th-order derivative is undefined. Therefore, one may prefer to have the minimum number of active knots, i.e., minimum number of jumps in the $(k-1)$th-order derivative of the B-spline curve. To do this end, one-norm minimization can be used which is considered as a method to result in solutions with only a few nonzero components, namely sparse solutions [41], [42]. One can formulate this problem as a term in the objective function as follows:

$$J_3 = \sum_{i=0}^{h-1} \frac{|x_{i+1}^{(k-1)} - x_i^{(k-1)}| + |y_{i+1}^{(k-1)} - y_i^{(k-1)}| + |z_{i+1}^{(k-1)} - z_i^{(k-1)}|}{T_s} \qquad (13)$$

*4) Reconnaissance zones:* In various missions [5], [43], [44], it is required for the UAV to pass through some preflight or in-flight midcourse waypoints, which can be considered as temporary targets. The preflight defined waypoints are specified according to the particular mission, which includes reconnaissance or implementation of some tasks [45], [46]. On the other hand in several missions there is no prior knowledge of the waypoints and they should be computed as an in-flight process and thus not programmable beforehand, or provided by the ground control stations [3].

In this study, the waypoints are added to the problem in the objective function in lieu of the constraints. Each waypoint is a temporary target, and thus can be considered instead of the final state in the obejctive function as follows:

$$J_1 = w_1 \sum_{i=0}^{h-1} ||\mathbf{x}_{i+1} - \mathbf{x}_i||_2 + w_2||\mathbf{x}_h - \mathbf{x}_{w_j}||_2 \qquad (14)$$

in which $\mathbf{x}_{w_j}$ is the position vector of the $j$th waypoint. Once the UAV reaches the $j$th waypoint, the objective function is modified to include the next waypoint as the final state. This procedure continues until the UAV passes through all the waypoints and finally reaches the endpoint. In some studies [47], it is favorable to define the reconnaissance zones as spherical regions of radius $R_z$ centered at the waypoints locations, and thus it is sufficient for the UAV to enter these zones and not necessarily reaching the waypoints. In order to incorporate the latter case, a stop condition is added to the problem in which the receding horizon procedure is stopped for each waypoint when,

$$||\mathbf{x}_1 - \mathbf{x}_{w_j}||_2 \leq R_z \qquad (15)$$

which ensures that the UAV enters the reconnaissance zones.

### B. Problem Constraints

This subsection touches on the formulation of the constraints involved in the related problem. These constraints incorporate dynamic (velocity, acceleration, and radius of curvature) ones, initial conditions, stationary obstacle avoidance, continuity conditions, terrain constraints, radar detection, flight-prohibited zones, and UAVs collision avoidance.

*1) Dynamic constraints:* The dynamic constraints considered in this paper are the maximum allowable velocity and acceleration, as well as the minimum allowable radius of curvature, which are mainly constrained by the maximum thrust of the motors [48]. Also, the minimum value for the velocity of a UAV is limited by some factors such as the stall effect and angle of attack.

The first derivative of a B-spline curve is obtained by [49]:

$$\dot{\mathbf{x}}(t) = \dot{\mathbf{r}}(t) = \sum_{l=0}^{n} \mathbf{q}_l \dot{N}_{l,k}(t) \qquad (16)$$

where,

$$\dot{N}_{l,k}(t) = \frac{k-1}{t_{l+k-1} - t_l} N_{l,k-1}(t) - \frac{k-1}{t_{l+k} - t_{l+1}} N_{l+1,k-1}(t) \qquad (17)$$

Thus the maximum allowable velocity constraint can be formulated as the following convex inequality:

$$||\dot{\mathbf{x}}_i||_2 \leq v_{\max}, \ i = 1, ..., h \tag{18}$$

where $v_{\max}$ is maximum feasible velocity for the UAV.

In order to constrain the velocity of the UAV to take values larger than a minimum allowable velocity, one can make the following suggestion:

$$||\dot{\mathbf{x}}_i||_2 \geq v_{\min}, \ i = 1, ..., h \tag{19}$$

Since the above inequalities are non-convex, the following more conservative substitutions are proposed based on the aforementioned disjunctive programming:

$$\begin{aligned}
&\dot{x}_i \geq v_{\min} \vee \dot{x}_i \leq -v_{\min} \vee \\
&\dot{y}_i \geq v_{\min} \vee \dot{y}_i \leq -v_{\min} \vee \\
&\dot{z}_i \geq v_{\min} \vee \dot{z}_i \leq -v_{\min}, \ i = 1, ..., h
\end{aligned} \tag{20}$$

where $\dot{\mathbf{x}}_i = [\dot{x}_i, \dot{y}_i, \dot{z}_i]^T$. Using the big-M method one can modify the above to mixed integer linear constraints for each $i = 1, ..., h$ as follows:

$$\begin{aligned}
-\dot{x}_i + v_{\min} &\leq (1 - b_{i1})M \\
\dot{x}_i + v_{\min} &\leq (1 - b_{i2})M \\
-\dot{y}_i + v_{\min} &\leq (1 - b_{i3})M \\
\dot{y}_i + v_{\min} &\leq (1 - b_{i4})M \\
-\dot{z}_i + v_{\min} &\leq (1 - b_{i5})M \\
\dot{z}_i + v_{\min} &\leq (1 - b_{i6})M \\
\sum_{j=1}^{6} b_{ij} &= 1
\end{aligned} \tag{21}$$

where $M$ is a sufficiently large constant and $b_{ij}$ is a binary variable. In some studies [50], it is desired to constrain the minimum and maximum velocities of the UAV in the horizontal plane ($xoy$) as well as that of the climb and descent (in $z$ direction). It is noteworthy that one can readily modify the proposed approach in this paper to the latter problem, using the following constraints:

$$\begin{aligned}
v_{h_{\min}} &\leq ||(\dot{x}_i, \dot{y}_i)||_2 \leq v_{h_{\max}} \\
|\dot{z}_i| &\leq v_{z_{\max}}
\end{aligned} \tag{22}$$

where $v_{h_{\min}}$, $v_{h_{\max}}$ and $v_{z_{\max}}$ are the minimum and maximum velocities in the horizontal plane and the maximum value of climb and descent rate, respectively.

Acceleration constraint can be implemented using the second derivative of B-spline functions, namely,

$$||\ddot{\mathbf{x}}_i||_2 \leq a_{\max}, \ i = 1, ..., h \tag{23}$$

where $a_{\max}$ is maximum possible acceleration for the UAV, and $\ddot{\mathbf{x}}(t)$ can be obtained from Eqs. (16) and (17) by induction.

In the following, we show that by constraining the magnitudes of the velocity and acceleration vectors, one can limit the maximum values of the curvature. To do this end, firstly, the curvature ($\kappa$) of a three dimensional curve, $\mathbf{r}(t)$, is defined by:

$$\kappa = \frac{1}{\rho} = ||\frac{d\mathbf{T}}{ds}||_2 = \frac{||\dot{\mathbf{T}}||_2}{\dot{s}} \tag{24}$$

where $\rho$ is the radius of curvature, $\mathbf{T}(t)$ is the tangent vector, and $s(t)$ is the arc length. One can readily prove that [51],

$$\kappa = \frac{||\dot{\mathbf{r}} \times \ddot{\mathbf{r}}||_2}{||\dot{\mathbf{r}}||_2^3} \tag{25}$$

From the above equation, it reveals that the curvature will have the greatest value if the acceleration vector is perpendicular to the velocity one, and therefore limited by,

$$\kappa \leq \left(\frac{||\ddot{\mathbf{r}}||_2}{||\dot{\mathbf{r}}||_2^2} = \frac{a}{v^2}\right) \leq \frac{a_{\max}}{v_{\min}^2} \tag{26}$$

where $v$ and $a$ are the velocity and acceleration values at each time. As a result we have,

$$\kappa \leq \frac{a_{\max}}{v_{\min}^2} \tag{27}$$

or,

$$\rho \geq \frac{v_{\min}^2}{a_{\max}} \tag{28}$$

Therefore, by constraining the magnitudes of the velocity and acceleration vectors, it is possible to limit the maximum value of the curvature or, equivalently, the minimum value of radius of curvature. Since it is favorable for the UAV to reach the target in the minimum time, it should be remained at or near the maximum velocity throughout its whole maneuver [52]–[54]. In this case, one can conclude that,

$$\kappa \approx \frac{a}{v^2} \tag{29}$$

and thus,

$$\kappa_{\max} \approx \frac{a_{\max}}{v_{\max}^2} \tag{30}$$

*2) Initial and continuity conditions:* Assuming the starting state of the UAV as $\mathbf{x}_s$, one can formulate the initial conditions over the first look-ahead horizon as follows:

$$\begin{aligned}
\mathbf{x}_0 &= \mathbf{x}_s \\
\dot{\mathbf{x}}_0 &= \dot{\mathbf{x}}_s \\
\ddot{\mathbf{x}}_0 &= 0
\end{aligned} \tag{31}$$

where $\dot{\mathbf{x}}_s$ is the starting velocity (takeoff velocity) of the UAV. If the initial values for the heading angle ($\phi$) and flight-path angle ($\psi$) are given in lieu of the takeoff velocity, one can substitute the second constraint in Eq. (31) with the following one:

$$\dot{\mathbf{x}}_0 = [v_s \cos\psi_s \cos\phi_s, v_s \cos\psi_s \sin\phi_s, v_s \sin\psi_s]^T \tag{32}$$

where $\psi_s$, $\phi_s$, and $v_s$ are the initial values for heading angle, flight-path angle, and magnitude of velocity, respectively.

It should be noted that the position, velocity, and acceleration vectors at the end of first time step of each horizon are the initial conditions for the subsequent one, i.e.,

$$\begin{aligned}
\mathbf{x}_0|_{t^+} &= \mathbf{x}_1|_{t^-} \\
\dot{\mathbf{x}}_0|_{t^+} &= \dot{\mathbf{x}}_1|_{t^-} \\
\ddot{\mathbf{x}}_0|_{t^+} &= \ddot{\mathbf{x}}_1|_{t^-}
\end{aligned} \tag{33}$$

By these conditions, we guarantee the parametric continuity of the obtained trajectory up to the second order.

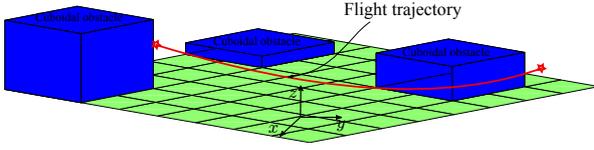

Fig. 2. Schematic showing the cuboidal obstacles and the fight trajectory.

*3) Obstacle avoidance:* In this study, statical obstacles with various shapes, which can be approximated with covering convex polyhedra, can be taken into consideration. Therefore, the obstacle avoidance formulation in this paper is only proposed for polyhedral obstacles. In Fig. 2, different rectangular cuboids are considered as the obstacles presented in the environment, which the UAV should not collide with.

In this paper, the convex polyhedron, $\mathcal{P}$, is defined as the intersection of a finite number of halfspaces, namely,

$$\mathcal{P} = \{\boldsymbol{\xi} | \mathbf{A}\boldsymbol{\xi} \preceq \mathbf{c}\}. \tag{34}$$

where,

$$\mathbf{A} = \begin{bmatrix} \mathbf{a}_1^T \\ \vdots \\ \mathbf{a}_p^T \end{bmatrix}, \quad \mathbf{c} = \begin{bmatrix} c_1 \\ \vdots \\ c_p \end{bmatrix} \tag{35}$$

in which $\mathbf{a}_i$ and $c_i$ specify the $i$th halfspace, and the symbol $\preceq$ denotes the componentwise inequality in $\mathbf{R}^p$, i.e., the real coordinate space of $p$ dimensions. The point $\mathbf{x}_i$ lies outside the polyhedron in case at least one of the $p$ inequalities in $\mathbf{A}\mathbf{x}_i \succeq \mathbf{c}$ is satisfied, which results in a disjunctive constraint. Using big-M method, one can modify this constraint to mixed integer linear ones, namely,

$$\begin{aligned} \mathbf{A}\mathbf{x}_i &\geq \mathbf{c} + \eta\mathbf{1} + (\mathbf{b}_i - \mathbf{1})M \\ \sum_{j=1}^{p} b_{ij} &= 1 \end{aligned} \tag{36}$$

where $\mathbf{b}_i = [b_{i1}, ..., b_{ip}]^T$ is the vector of binary variables. Moreover, $\eta \geq 0$ is a safety margin, which causes that a minimum distance be preserved from the boundaries of the obstacle.

By the proposed formulation, we can implement the obstacle-avoidance constraints. These obstacles can be those recognized beforehand during the reconnaissance operations, and provided by the ground control stations, or unknown ones detected by the sensors of the UAV.

It should be noted that in case of moving obstacles, $\mathbf{A}$ and $\mathbf{c}$ change with respect to the time. In other words, the obstacle polyhedron should be represented by the following set:

$$\mathcal{P}(t) = \{\boldsymbol{\xi} | \mathbf{A}(t).\boldsymbol{\xi} \preceq \mathbf{c}(t)\} \tag{37}$$

Tehrefore, having the instantaneous formulation of the obstacle polyhedron, one can implement the dynamic obstacle avoidance constraints.

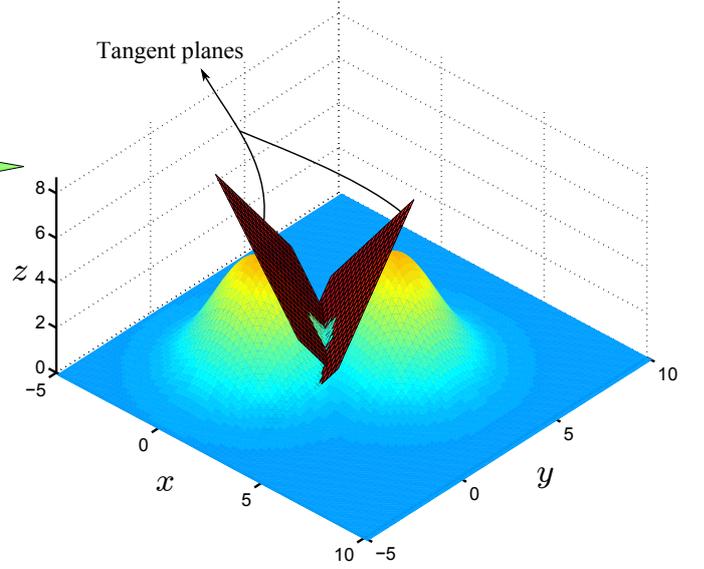

Fig. 3. Terrain collision avoidance using tangent planes strategy.

*4) Terrain constraints:* As mentioned previously, the feasible trajectory cannot go through the terrain, namely, there should not be any contact point with the terrain. To do this end, the tangent plane strategy is touched upon in the following.

This approach aims to find one or more planes which have resembling local features with those of the terrain using the data provided by the UAV sensors. For this purpose, first, the nearest point of the terrain to the UAV, considered as the most hazardous, is obtained, and then the tangent plane at this point is computed via the normal vector, which can be calculated by the numerical gradient method. This plane divides $\mathbf{R}^3$ into two halfspaces, one of which is considered as the safe region. To avoid collision with the terrain, it should be guarateed that the UAV stays within the aforementioned safe region during each horizon. Assuming that the safe region is represented by the following set:

$$\{\mathbf{x} \mid \mathbf{a}^T\mathbf{x} \leq c\} \tag{38}$$

then, the terrain constraint is obtained as:

$$\mathbf{a}^T\mathbf{x}_i - c + \kappa \leq 0, \ i = 1, ..., h \tag{39}$$

where $\kappa \geq 0$ is a safety margin, which causes that a minimum distance be preserved from the terrain.

It should be noted that it is possible to select more than one tangent plane. It is suggested that first a minimum distance between the terrain and the UAV, namely $d_{\min}$, is chosen, and then for each detected point of the terrain from which the distance of the UAV is less than $d_{\min}$, the tangent plane at that point is computed. For instance, as illustrated in Fig. 3, if the UAV should pass through two adjacent hillsides, the tangent planes at the two nearest points of the terrain to the UAV are obtained, and then it is restricted to move between these planes.

*5) Radar detection and flight-prohibited zones:* In addition to unknown radars that for escaping their detection the UAV should fly at low altitudes, there can be other particular radars in the environment, such as those in air defense units, that the UAV should avoid being detected by them, thereby the UAV and its mission will be kept safe. Furthermore, there may exist some special zones considered as high-risk zones that the UAV must not enter them.

Each of the above regions is defined as a rectangular zone with the external limits, i.e., $x_l$, $x_u$, $y_l$, and $y_u$. For the purpose of limiting the UAV to stay outside this region one can apply the disjunctive programming in addition to the big-M method in order to obtain the following required constraints for $i = 1, ..., h$:

$$\begin{aligned} x_i - x_l + \gamma &\leq (1 - b_{i1})M \\ x_u - x_i + \gamma &\leq (1 - b_{i2})M \\ y_i - y_l + \gamma &\leq (1 - b_{i3})M \\ y_u - y_i + \gamma &\leq (1 - b_{i4})M \\ \sum_{j=1}^{6} b_{ij} &= 1 \end{aligned} \quad (40)$$

where $M$ and $b_{ij}$ are, as before, a sufficiently large constant and a binary variable, respectively. In the above constraints, $\gamma \geq 0$ is a safety margin, which causes that a minimum distance be preserved from the boundaries of the no flying zone.

*6) UAVs collision avoidance:* The last constraint involved in this study is a cooperation one which becomes important in generating trajectories for multiple UAVs, in which it is should be checked whether two UAVs are getting too close in order to avoid collision while they are following their respective trajectories. This constraint can be formulated for any two UAVs with position vectors as $\mathbf{x}(t)$ and $\mathbf{x}'(t)$ using a minimum allowable distance, $d_{all}$, between them as follows:

$$||\mathbf{x}_i - \mathbf{x}'_i||_2 \geq d_{all}, \ i = 1, ..., h \quad (41)$$

The above are non-convex constraints and thus, a relaxation strategy is indispensably required in order to convexify them. In this regard, the disjunctive programming is employed to play the role of a convexifier by applying the follwoing mixed integer linear constraints instead of the aforementioned non-convex ones:

$$\begin{aligned} \mathbf{x}_i - \mathbf{x}'_i + d_{all} &\leq (1 - b_{i1})M \\ \mathbf{x}'_i - \mathbf{x}_i + d_{all} &\leq (1 - b_{i2})M \\ \mathbf{y}_i - \mathbf{y}'_i + d_{all} &\leq (1 - b_{i3})M \\ \mathbf{y}'_i - \mathbf{y}_i + d_{all} &\leq (1 - b_{i4})M \\ \mathbf{z}_i - \mathbf{z}'_i + d_{all} &\leq (1 - b_{i5})M \\ \mathbf{z}'_i - \mathbf{z}_i + d_{all} &\leq (1 - b_{i6})M \\ \sum_{j=1}^{6} b_{ij} &= 1 \end{aligned} \quad (42)$$

where $M$ is a sufficiently large constant and $b_{ij}$ is a binary variable. It should be noted that if the maximum distance between the UAVs is favorable, one can add $d_{all}$ as a variable to the objective funtion which should be maximized or equivalently minimizing $-d_{all}$.

In some studies it is required that the UAVs stay within a communication range of one another which depends on the relative distance between them [55]. Assuming that the communication range is a spherical region with radius $R_c$ centered at the location of each UAV, one can add the following convex constraints to the problem, which restrict the maximum allowable distance between the UAVs:

$$||\mathbf{x}_i - \mathbf{x}'_i||_2 \leq R_c, \ i = 1, ..., h \quad (43)$$

TABLE I
PHYSICAL CHARACTERISTICS OF THE UNDER-STUDY UAV.

| | |
|---|---|
| $v_{\max}$ | 60 (m/s) |
| $v_{h_{\min}}$ | 30 (m/s) |
| $v_{z_{\max}}$ | 6 (m/s) |
| $a_{\max}$ | 1.5g ($g = 9.81$ (m/s$^2$)) |
| $R$ | 1000 (m) |
| $\theta_1$ | 60° |
| $\theta_2$ | 45° |

## VI. SIMULATION RESULTS

To verify the performance of the proposed flight trajectory planning methodology, different numerical simulations are provided in this section. The under-study UAV has the physical characteristics shown in Table I. It should be noted that all of these examples have been implemented and solved via the solver Gurobi [39] for $t_h = 10$ (s) and $T_s = 1$ (s). In the first scenario, the optimal trajectory is obtained for different objective functions proposed in Section 5 while the coordinates of the starting and final points are specified as follows:

$$\mathbf{x}_s = [-2000, -2000, 10]^T, \ \mathbf{x}_f = [715, 1730, 4]^T$$

The problem is solved for four objective functions, namely $J = J_1$, $J = J_1 + J_3$, $J = J_1 + J_2$, and $J = J_1 + J_2 + J_3$, for which the optimal trajectories are represented in Fig. VI by T$_1$, T$_2$, T$_3$, and T$_4$, respectively. Furthermore, to verify that the dynamic constraints have been satisfied, the velocity and acceleration graphs for each trajectory are represented in Figs. 5, 6, 7, and 8. The four illustrated graphs for each trajectory consist of the total velocity ($v$), velocity in the horizontal plane ($v_h$), rate of climb and descent ($v_z$), and the acceleration ($a$) of the UAV with respect to time. In generating these trajectories, for preventing the UAV from colliding with the terrain, the tangent plane strategy is taken into consideration while the safety margin is selected as $\kappa = 1$ (m). The lengths of the four optimal trajectories obtained in this problem have been represented in Table II.

It is obvious from the results that T$_1$ and T$_2$ have shorter lengths than those of T$_3$ and T$_4$, which is due to the fact that in the latter the UAV should fly at low altitudes. In other words, the accumulated altitude of the UAV during the time horizon is minimized in order not to be detected by unknown radars. However, as mentioned previously, flying at low altitude can make the UAV travel a longer trajectory which has the drawback of increasing the fuel consumption. In this case, flying at low altitude caused an almost 400 (m) increase in the length of the trajectory that the UAV needs to fly before reaching its target which prolongs the flight duration as observable in the velocity and acceleration graphs.

TABLE II
LENGTHS OF THE FOUR OPTIMAL TRAJECTORIES.

| Trajectory | Length (m) |
| --- | --- |
| $T_1$ | 4663.8 |
| $T_2$ | 4663.3 |
| $T_3$ | 5086.5 |
| $T_4$ | 5037.8 |

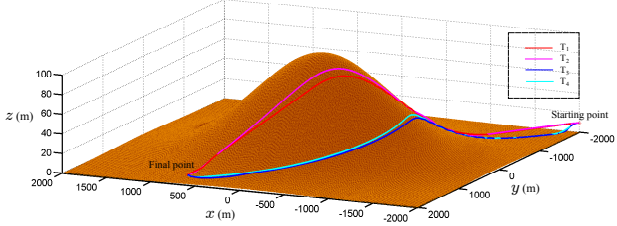

Fig. 4. Optimal trajectories obtained for different objective functions.

It should be noted that all the generated motion trajectories are $C^2$ continuous. Moreover, as stated before, it is favorable to have the minimum number of active knots, i.e., minimum number of jumps in the third-order derivatives of the trajectories using one-norm minimization. By comparing the results illustrated in Figs. 5 and 7 with those depicted in Figs. 6 and 8, it reveals that the number of fluctuations in velocity graphs and that of breaking points in acceleration graphs have been reduced tremendously due to the one-norm minimization which is added to the objective function through the term $J_3$.

In the second scenario, the coordinates of the starting and final points are specified as follows:

$$\mathbf{x}_s = [3000, -2000, 10]^T, \quad \mathbf{x}_f = [-2000, 3000, 1]^T$$

Moreover, it is assumed that there is a cuboidal obstacle as well as a flight-prohibited zone in the enviornment which should be avioded by the UAV and their geometric characteristics are specified in Table III. This problem is formulated based on the approach presented in Section 5, and then solved for which the result is represented in Fig. 9. It is noteworthy that the objective function is selected as $J = J_1 + J_2 + J_3$, and the following safety margins are taken into account:

$$\kappa = 1 \text{ (m)}, \quad \eta = 5 \text{ (m)}, \quad \gamma = 20 \text{ (m)}$$

It can be observed in Fig. 9 that for the obtained trajectory the UAV will fly over the cuboidal obstacle, and then it will turn the flight-prohibited zone since it is forbiddden to pass over it. Due to space limitation the velocity and acceleration graphs are not represented for the second and third scenarios.

The next scenario aims at testing the proposed method for multiple UAVs. In this regard, it is assumed that in this mission there are two UAVs and the starting and final points for each are specified as follows:

$$\mathbf{x}_{s_1} = [-200, -2000, 10]^T, \quad \mathbf{x}_{f_1} = [200, 2000, 4]^T$$
$$\mathbf{x}_{s_2} = [200, -2000, 10]^T, \quad \mathbf{x}_{f_2} = [-200, 2000, 4]^T$$

It is considered a symmetry between the starting and final points of the two UAVs in order to show how they can escape collision with each other. The objective function is selected as $J = J_1 + J_2 + J_3$ and the safety margin is chosen as $\kappa = 1$ (m). This problem is formulated and then solved according to the approach presented in Section 5, for which the result is represented in Fig. 10. It is obvious in Fig. 10 that UAV 2 changes its path by climbing and

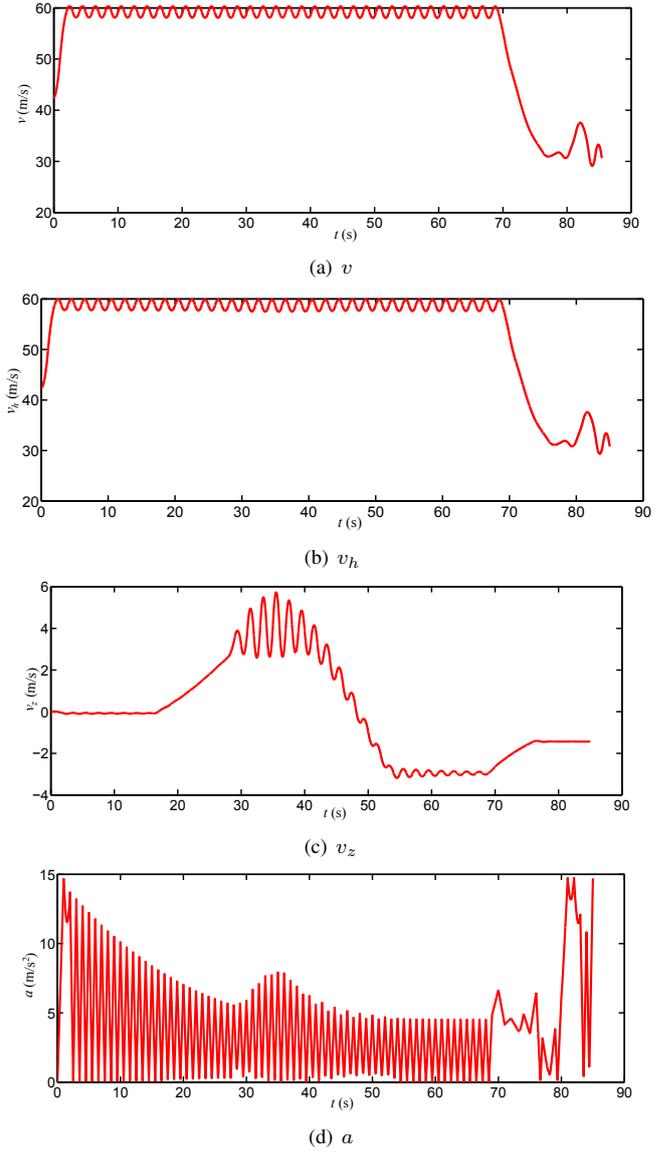

(a) $v$

(b) $v_h$

(c) $v_z$

(d) $a$

Fig. 5. Velocity and acceleration graphs for $T_1$.

TABLE III
GEOMETRIC CHARACTERISTICS OF THE CUBOIDAL OBSTACLE AND FLIGHT-PROHIBITED ZONE.

(a) Cuboidal obstacle

| $x_l$ | 1600 |
| --- | --- |
| $y_l$ | -1400 |
| $z_l$ | 0 |
| $x_u$ | 2400 |
| $y_u$ | -600 |
| $z_u$ | 50 |

(b) Flight-prohibited zone

| $x_l$ | -1500 |
| --- | --- |
| $y_l$ | 2000 |
| $x_u$ | 2400 |
| $y_u$ | -600 |

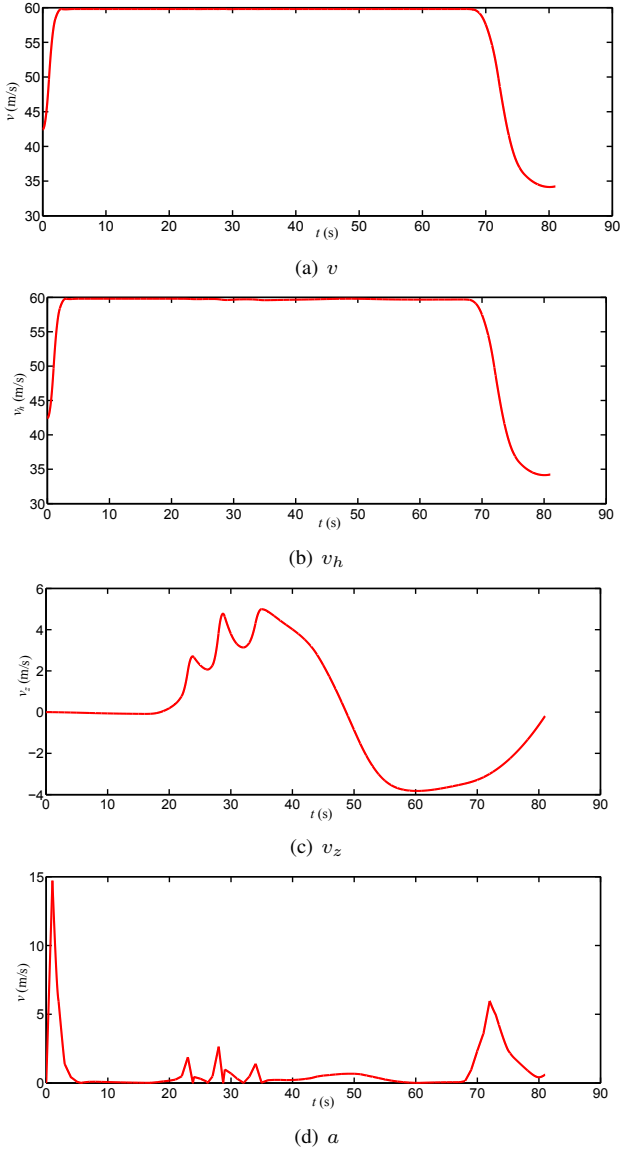

Fig. 6. Velocity and acceleration graphs for T$_2$.

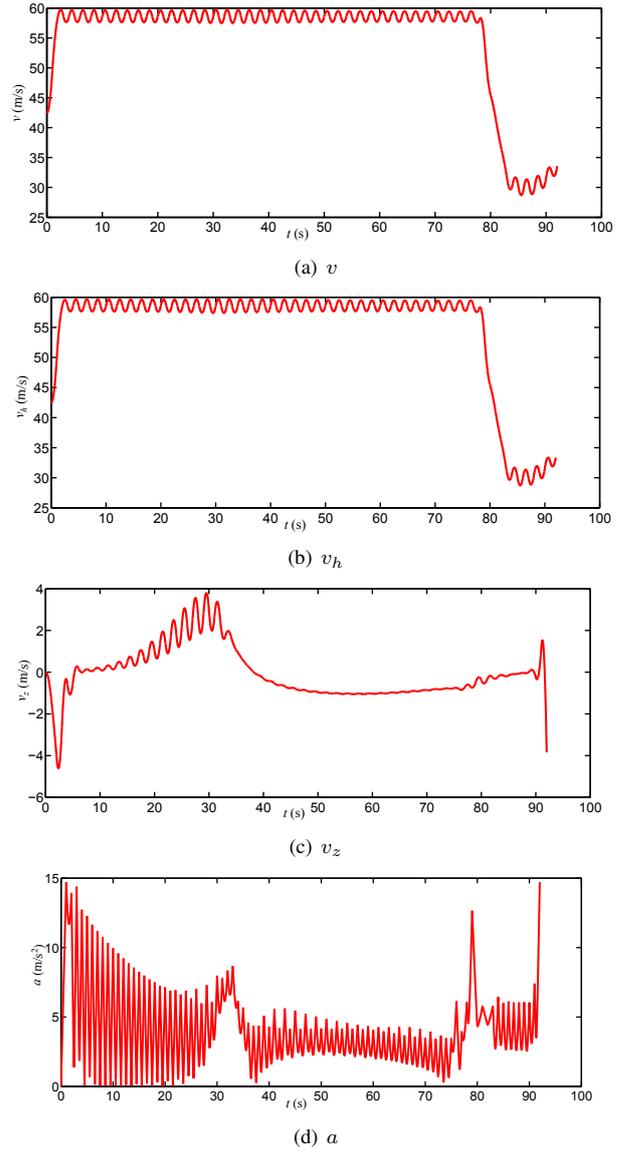

Fig. 7. Velocity and acceleration graphs for T$_3$.

then descending in order to avoid collision with UAV 1. It should be noted that the minimum and maximum allowable distances between the two UAVs are assumed as $d_{all} = 50$ (m) and $R_c = 400$ (m), respectively. Figure 11 depicts the instantaneous distance between the UAVs which confirms the aforementioned allowable distances between them. Finally, in order to show how the collision avoidance constraint can contribute to this problem, in Fig. 12 the optimal trajectories are displayed without considering collision avoidance constraint which reveals a collision between the two UAVs.

In the last scenario, the coordinates of the starting and final points are specified as follows:

$$\mathbf{x}_s = [3000, -2000, 10]^T, \quad \mathbf{x}_f = [-2500, 3500, 1]^T$$

Moreover, it is assumed that there is a dynamic obstacle in the enviornment which should be avoided by the UAV. This problem is formulated based on the approach presented in Section 5, and then solved for which the results are represented in Figs. 6 and 7. It is noteworthy that the objective function is selected as $J = J_1 + J_2 + J_3$.

It can be observed in Fig. 13 that for the obtained trajectory the UAV will fly over the dynamic obstacle in order to avoid collision. In addition, the problem is solved without considering the presence of the dynamic obstacle for which the result is illustrated in Fig. 14. One can observe that in this case the UAV collides with the obstacle.

It is noteworthy that the computational time for solving the optimization problem in each horizon for above scenarios is approximately in the interval of $[0.2 \text{ (s)}, 2 \text{ (s)}]$ using Matlab code on a PC equipped with an Intel(R) Core(TM) i5-2430M CPU @ 2.40GHz, and 4GB RAM.

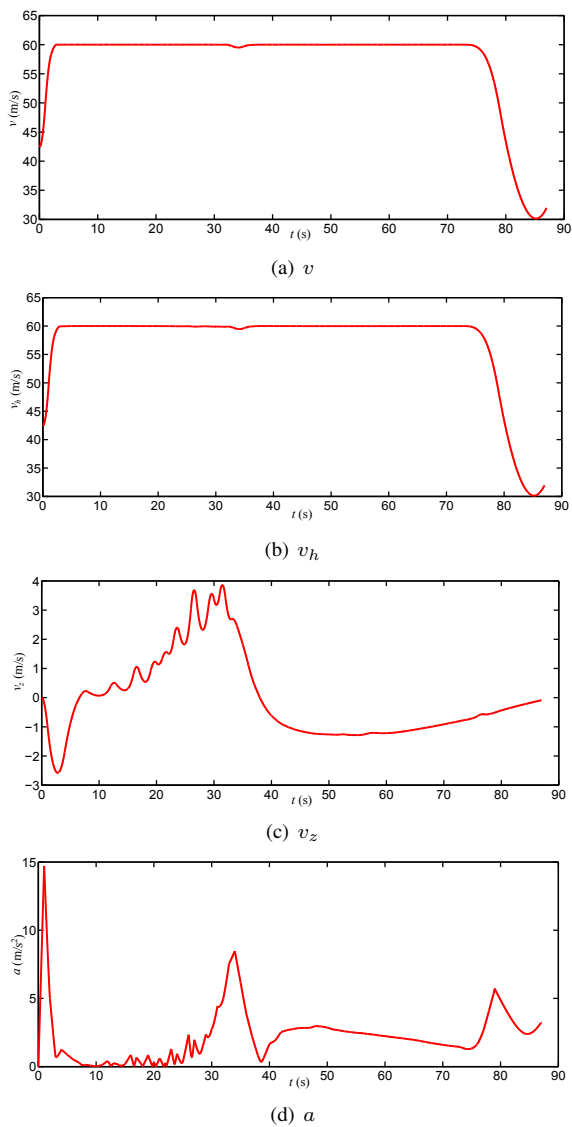

(a) $v$

(b) $v_h$

(c) $v_z$

(d) $a$

Fig. 8. Velocity and acceleration graphs for T$_4$.

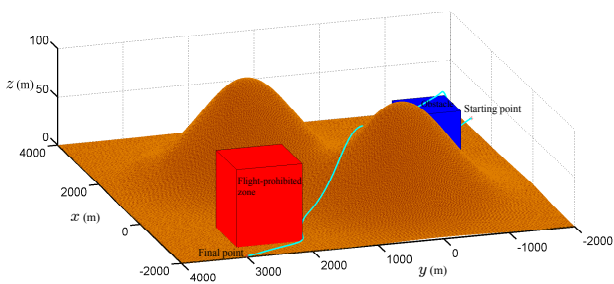

Fig. 9. Optimal trajectory obtained in the presence of obstacle and flight-prohibited zone.

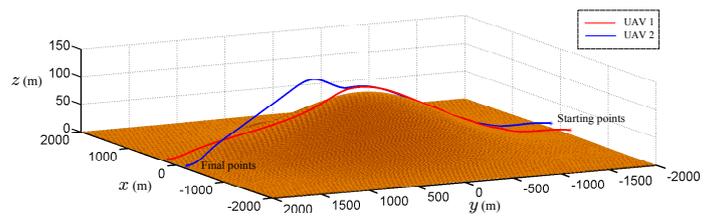

Fig. 10. Optimal trajectories obtained for two UAVs.

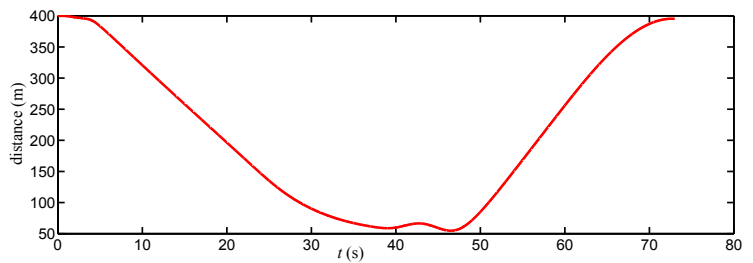

Fig. 11. Distance between the two UAVs.

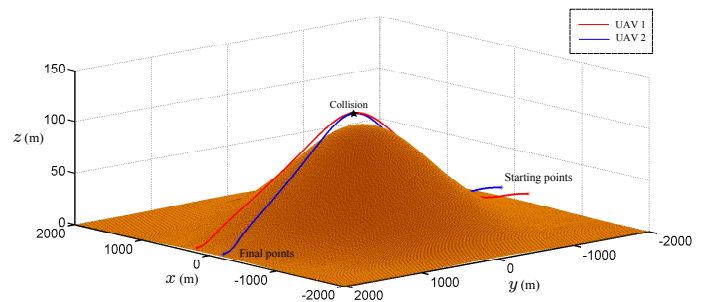

Fig. 12. Optimal trajectories obtained for two UAVs without considering collision avoidance constraint.

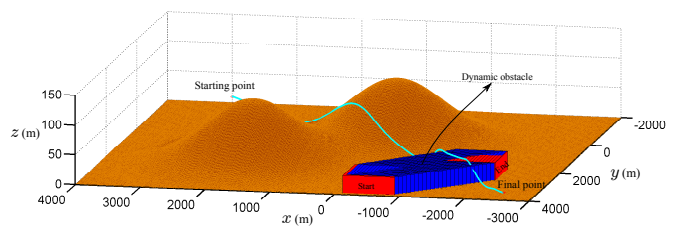

Fig. 13. Optimal trajectory obtained for the UAV in the presence of dynamic obstacle.

## VII. CONCLUSION

In this paper the problem of trajectory planning for multiple UAVs was investigated in the presence of static obstacles, three-dimensional terrains, radar detection, and flight-prohibited zones using a combination of disjunctive optimization and receding horizon concept. Moreover, a strategy is proposed in order to constrain the UAV to pass through some preflight or in-flight midcourse waypoints. The fourth-order B-spline curves are used in this paper and it is ensured that the generated trajectories have at least $C^2$ continuity. Furthermore, the multi-objective function incorporates different terms such as trajectory length, flight altitude, and number of active knots which increases the smoothness of the generated curves. Thereafter, the results for different test scenarios obtained via the method presented in this paper. Using the concepts of dis-crete optimization, convex optimization, and receding horizon lead to a reduction in the computational time, which made the algorithm appropriate for online applications. Modication of the proposed algorithm in the presence of wind can be considered as the ongoing study.

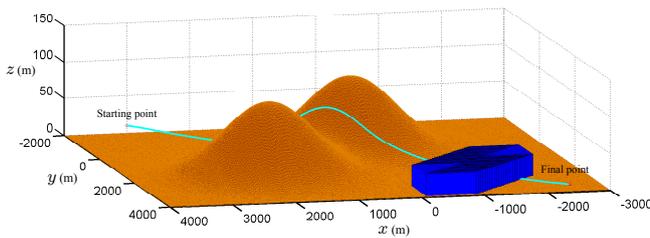

Fig. 14. Optimal trajectory obtained for the UAV without considering the constraint of dynamic obstacle.

## REFERENCES

[1] J. F. Keane and S. S. Carr, "A brief history of early unmanned aircraft," *Johns Hopkins APL Technical Digest*, vol. 32, no. 3, pp. 558–571, 2013.
[2] K. Dalamagkidis, "Aviation history and unmanned flight," in *Handbook of Unmanned Aerial Vehicles*. Springer, 2014, pp. 57–81.
[3] A. Babaei and M. Mortazavi, "Three-dimensional curvature-constrained trajectory planning based on in-flight waypoints," *Journal of aircraft*, vol. 47, no. 4, pp. 1391–1398, 2010.
[4] A. Ollero and I. Maza, *Multiple heterogeneous unmanned aerial vehicles*. Springer Publishing Company, Incorporated, 2007.
[5] F. Nex and F. Remondino, "UAV for 3d mapping applications: a review," *Applied Geomatics*, vol. 6, no. 1, pp. 1–15, 2014.
[6] P. C. Merrell, D.-J. Lee, and R. W. Beard, "Obstacle avoidance for unmanned air vehicles using optical flow probability distributions," *Mobile Robots XVII*, vol. 5609, no. 1, pp. 13–22, 2004.
[7] I. Nikolos, N. Tsourveloudis, and K. Valavanis, "Evolutionary algorithm based path planning for multiple UAV cooperation," in *Advances in Unmanned Aerial Vehicles*. Springer, 2007, pp. 309–340.
[8] R. J. Szczerba, "Threat netting for real-time, intelligent route planners," in *Information, Decision and Control, 1999. IDC 99. Proceedings. 1999*. IEEE, 1999, pp. 377–382.
[9] E. Besada-Portas, L. de la Torre, J. M. de la Cruz, and B. de Andrés-Toro, "Evolutionary trajectory planner for multiple UAVs in realistic scenarios," *Robotics, IEEE Transactions on*, vol. 26, no. 4, pp. 619–634, 2010.
[10] H. Moradian and G. Vossoughi, "Robust velocity control for an a-shaped micro-robot with stick-slip locomotion," *Proceedings of the Institution of Mechanical Engineers, Part C: Journal of Mechanical Engineering Science*, vol. 230, no. 14, pp. 2413–2426, 2016.
[11] J. J. Ruz, O. Arévalo, J. M. de la Cruz, and G. Pajares, "Using milp for UAVs trajectory optimization under radar detection risk," in *Emerging Technologies and Factory Automation, 2006. ETFA'06. IEEE Conference on*. IEEE, 2006, pp. 957–960.
[12] A. L. Jennings, R. Ordónez, and N. Ceccarelli, "Dynamic programming applied to UAV way point path planning in wind," in *Computer-Aided Control Systems, 2008. CACSD 2008. IEEE International Conference on*. IEEE, 2008, pp. 215–220.
[13] J. M. de la Cruz, E. Besada-Portas, L. Torre-Cubillo, B. Andres-Toro, and J. A. Lopez-Orozco, "Evolutionary path planner for UAVs in realistic environments," in *Proceedings of the 10th annual conference on Genetic and evolutionary computation*. ACM, 2008, pp. 1477–1484.
[14] C. Zheng, L. Li, F. Xu, F. Sun, and M. Ding, "Evolutionary route planner for unmanned air vehicles," *Robotics, IEEE Transactions on*, vol. 21, no. 4, pp. 609–620, 2005.
[15] M. Pelosi, C. Kopp, and M. Brown, "Range-limited UAV trajectory using terrain masking under radar detection risk," *Applied Artificial Intelligence*, vol. 26, no. 8, pp. 743–759, 2012.
[16] T. Sharma, P. Williams, C. Bil, and A. Eberhard, "Optimal three dimensional aircraft terrain following and collision avoidance," *ANZIAM Journal*, vol. 47, pp. 695–711, 2007.
[17] Y. Chen, W. Zu, G. Fan, and H. Chang, "Unmanned aircraft vehicle path planning based on svm algorithm," in *Foundations and Practical Applications of Cognitive Systems and Information Processing*. Springer, 2014, pp. 705–714.
[18] S. Malaek and A. Abbasi, "Near-optimal terrain collision avoidance trajectories using elevation maps," *Aerospace and Electronic Systems, IEEE Transactions on*, vol. 47, no. 4, pp. 2490–2501, 2011.
[19] I. J. Schönberg, "Contributions to the problem of approximation of equidistant data by analytic functions," *Quart. Appl. Math*, vol. 4, no. 2, pp. 45–99, 1946.
[20] M. G. Cox, "The numerical evaluation of b-splines," *IMA Journal of Applied Mathematics*, vol. 10, no. 2, pp. 134–149, 1972.
[21] C. De Boor, "On calculating with b-splines," *Journal of Approximation Theory*, vol. 6, no. 1, pp. 50–62, 1972.
[22] ——, "A practical guide to splines," *Mathematics of Computation*, 1978.
[23] G. Farin, *Curves and surfaces for computer-aided geometric design: a practical guide*. Elsevier, 2014.
[24] I. W. Böhm, "Cubic b-spline curves and surfaces in computer aided geometric design," *Computing*, vol. 19, no. 1, pp. 29–34, 1977.
[25] L. Schumaker, *Spline functions: basic theory*. Cambridge University Press, 2007.
[26] Y. Kuwata, "Real-time trajectory design for unmanned aerial vehicles using receding horizon control," Ph.D. dissertation, Massachusetts Institute of Technology, 2003.
[27] D. Salomon, *Curves and surfaces for computer graphics*. Springer Science & Business Media, 2007.
[28] N. M. Patrikalakis and T. Maekawa, *Shape interrogation for computer aided design and manufacturing*. Springer Science & Business Media, 2009.
[29] ——, "Intersection problems," in *Handbook of computer aided geometric design*. Elsevier, 2002, pp. 623–649.
[30] V. Rovenski, "Modeling of curves and surfaces with matlab," ser. Springer Undergraduate Texts in Mathematics and Technology. Springer New York, 2010.
[31] S. P. Boyd and L. Vandenberghe, *Convex optimization*. Cambridge university press, 2004.
[32] T. Schouwenaars, J. How, and E. Feron, "Receding horizon path planning with implicit safety guarantees," in *American Control Conference, 2004. Proceedings of the 2004*, vol. 6. IEEE, 2004, pp. 5576–5581.
[33] T. Schouwenaars, É. Féron, and J. How, "Safe receding horizon path planning for autonomous vehicles," in *PROCEEDINGS OF THE ANNUAL ALLERTON CONFERENCE ON COMMUNICATION CONTROL AND COMPUTING*, vol. 40, no. 1. The University; 1998, 2002, pp. 295–304.
[34] L. Blackmore and B. Williams, "Optimal manipulator path planning with obstacles using disjunctive programming," in *American Control Conference, 2006*. IEEE, 2006, pp. 3–pp.
[35] I. E. Grossmann, "Review of nonlinear mixed-integer and disjunctive programming techniques," *Optimization and Engineering*, vol. 3, no. 3, pp. 227–252, 2002.
[36] I. E. Grossmann and S. Lee, "Generalized convex disjunctive programming: Nonlinear convex hull relaxation," *Computational Optimization and Applications*, vol. 26, no. 1, pp. 83–100, 2003.
[37] N. W. Sawaya and I. E. Grossmann, "A cutting plane method for solving linear generalized disjunctive programming problems," *Computers & chemical engineering*, vol. 29, no. 9, pp. 1891–1913, 2005.
[38] I. I. Cplex, "12.2 users manual," 2010.
[39] Optimization, "Gurobi, Inc. Gurobi Optimizer 6.0."
[40] T. Jameson, "A fuel consumption algorithm for unmanned aircraft systems," DTIC Document, Tech. Rep., 2009.
[41] B. Demeulenaere, G. Pipeleers, J. De Caigny, J. Swevers, J. De Schutter, and L. Vandenberghe, "Optimal splines for rigid motion systems: a convex programming framework," *Journal of Mechanical Design*, vol. 131, no. 10, p. 101004, 2009.
[42] W. Van Loock, G. Pipeleers, J. De Schutter, and J. Swevers, "A convex optimization approach to curve fitting with b-splines," in *Proceedings of the 18th IFAC World Congress*, vol. 18, 2011, pp. 2290–2295.
[43] G. H. Elkaim, F. A. P. Lie, and D. Gebre-Egziabher, "Principles of guidance, navigation, and control of UAVs," *Handbook of Unmanned Aerial Vehicles*, pp. 347–380, 2015.
[44] D.-W. Gu, W. Kamal, and I. Postlethwaite, "A UAV waypoint generator," in *AIAA 1st Intelligent Systems Technical Conference, Chicago*, 2004.
[45] T. R. Jorris and R. G. Cobb, "Multiple method 2-d trajectory optimization satisfying waypoints and no-fly zone constraints," *Journal of Guidance, Control, and Dynamics*, vol. 31, no. 3, pp. 543–553, 2008.

[46] F. Baralli, L. Pollini, and M. Innocenti, "Waypoint-based fuzzy guidance for unmanned aircraft a new approach," in *AIAA Guidance, Navigation, and Control Conference and Exhibit*, 2002, pp. 5–8.

[47] J. L. Foo, J. Knutzon, V. Kalivarapu, J. Oliver, and E. Winer, "Path planning of unmanned aerial vehicles using b-splines and particle swarm optimization," *Journal of aerospace computing, Information, and communication*, vol. 6, no. 4, pp. 271–290, 2009.

[48] T. Puls, H. Winkelmann, S. Eilers, M. Brucke, and A. Hein, "Interaction of altitude control and waypoint navigation of a 4 rotor helicopter," in *Advances in Robotics Research*. Springer, 2009, pp. 287–298.

[49] J. Prochazkova, "Derivative of b-spline function," in *Proceedings of the 25th Conference on Geometry and Computer Graphics. Prague, Czech Republic*, 2005.

[50] Y. Kuwata, "Trajectory planning for unmanned vehicles using robust receding horizon control," Ph.D. dissertation, Massachusetts Institute of Technology, 2006.

[51] J. Casey, *Exploring curvature*. Vieweg Wiesbaden, 1996.

[52] A. Richards and J. P. How, "Aircraft trajectory planning with collision avoidance using mixed integer linear programming," in *American Control Conference, 2002. Proceedings of the 2002*, vol. 3. IEEE, 2002, pp. 1936–1941.

[53] Y. Kuwata, A. Richards, T. Schouwenaars, and J. P. How, "Distributed robust receding horizon control for multivehicle guidance," *Control Systems Technology, IEEE Transactions on*, vol. 15, no. 4, pp. 627–641, 2007.

[54] Y. Kuwata and J. How, "Three dimensional receding horizon control for UAVs," in *AIAA Guidance, Navigation, and Control Conference and Exhibit*, vol. 3, 2004, pp. 2100–2113.

[55] R. W. Beard and T. W. McLain, "Multiple UAV cooperative search under collision avoidance and limited range communication constraints," in *Decision and Control, 2003. Proceedings. 42nd IEEE Conference on*, vol. 1. IEEE, 2003, pp. 25–30.